\documentclass[12pt]{amsart}

\usepackage{amsmath,amsthm}
\usepackage{amssymb}
\usepackage{mathrsfs}
\usepackage{amsfonts}

\usepackage{enumerate}

\usepackage{graphicx}

\usepackage[T1]{fontenc}

\makeatletter
\@namedef{subjclassname@2010}{%
  \textup{2010} Mathematics Subject Classification}
\makeatother

\newtheorem{theorem}{Theorem}[section]
\newtheorem{corollary}[theorem]{Corollary}
\newtheorem{lemma}[theorem]{Lemma}
\newtheorem{proposition}[theorem]{Proposition}

\theoremstyle{definition}
\newtheorem{definition}{Definition}

\newtheorem{example}[theorem]{Example}

\newtheorem*{notation}{Notation}

\newcommand{\Z}{\mathbb{Z}}        
\newcommand{\R}{\mathbb{R}}        
\newcommand{\N}{\mathbb{N}}        
\newcommand{\Rn}{\mathbb{R}^d}     

\newcommand{\CN}{C^\infty(\Rn)}

\renewcommand{\DH}{\cD_H'(\Rn)}

\newcommand{\OC}{\cO_C'(\Rn)}
\newcommand{\Oc}{\cO_C'(\R)}
\newcommand{\OH}{\cO_H'(\Rn)}
\newcommand{\Oh}{\cO_H'(\R)}

\newcommand{\Di}{\cD(\Rn)}
\newcommand{\Dip}{\cD\,'(\Rn)}
\newcommand{\NZ}{\Rn_*}

\newcommand{\cD}{\mathscr{D}}

\newcommand{\cE}{\mathscr{E}}
\newcommand{\cS}{\mathscr{S}}
\newcommand{\cO}{\mathscr{O}}
\newcommand{\cA}{\mathscr{A}}

\newcommand{\cM}{\mathcal{M}}
\newcommand{\be}{{\bf 1}}
\newcommand{\supp}{\mathrm{supp}\, }
\renewcommand{\span}{\mathrm{span}}

\newcommand{\vp}{\varphi}
\newcommand{\eps}{\varepsilon}

\begin{document}


\baselineskip=17pt


\title{Hadamard operators on $\cD'(\Rn)$}

\author[D. Vogt]{Dietmar Vogt\\
Bergische Universit{\"a}t Wuppertal\\
FB Math.-Nat., Gauss-Str. 20,\\
D-42119 Wuppertal, Germany\\
E-mail: dvogt@math.uni-wuppertal.de}

\dedicatory{Dedicated to the memory of Pawe{\l} Doma\'nski}


\begin{abstract}
 We study continuous linear operators on $\Dip$ which admit all monomials as eigenvectors, that is, operators of Hadamard type. Such operators on $C^\infty(\mathbb{R}^d)$ and on the space $\cA(\mathbb{R}^d)$ of real analytic functions on $\mathbb{R}^d$ have been investigated by Doma\'nski, Langenbruch and the author. The situation in the present case, however, is quite different and also the characterization. An operator $L$ on $\Dip$ is of Hadamard type if there is a distribution $T$, the support of which has positive distance to all coordinate hyperplanes and which has a certain behaviour at infinity, such  that $L(S) = S\star T$ for all $S\in\Dip$. Here $(S\star T)\varphi = S_y(T_x\varphi(xy))$ for all $\varphi\in\Di$. To describe the behaviour at infinity we introduce a class $\mathscr{O}_H'(\mathbb{R}^d)$ of distributions defined by the same conditions like in the description of class $\mathscr{O}_C'(\mathbb{R}^d)$ of Laurent Schwartz, but derivatives replaced with Euler derivatives.
\end{abstract}



\maketitle

In the present note we study Hadamard operators on $\Dip$, that is, continuous linear operators on $\Dip$ which admit all monomials as eigenvectors and we give a complete characterization. Such operators on $\CN$ have been studied and characterized in \cite{Vhad, Vconv}, on $\cA(\R)$ in \cite{DLre, DLalg, DLhad} and on $\cA(\Rn)$ in \cite{DLV}. There you find also references to the long history of such problems. Since it can be shown that Hadamard operators commute with dilations our problem is, by duality, closely related to the study of continuous linear operators in $\Di$ which commute with dilations. They have the form $\vp\mapsto T_x\vp(xy)$ where $T$ is a distribution. The class $\DH$ of distributions $T$ such that $T_x\vp(xy)\in \Di$ for every $\vp\in\Di$ is studied. These are distributions with positive distance to the coordinate hyperplanes and certain behaviour at infinity, similar to the class $\cO_C'$ of L. Schwartz of rapidly decreasing distributions. We define a class $\OH$ of distributions by the same conditions like in the description of class $\cO_C'$ in \cite{LS}, but derivatives replaced with Euler derivatives. We denote by $\cM(\Rn)$ the class of Hadamard operators in $\Dip$ and for $x\in\NZ=(\R\setminus\{0\})^d$ we set $\sigma(x)=\prod_j\frac{x_j}{|x_j|}$, that is, the signum of $x$ which, of course, is constant on each `{quadrant}'. Then our main result is (see Corollary \ref{c1}, Theorem \ref{t11} and Theorem \ref{t7}):

{\bf Main Theorem} \it $L\in\cM(\Rn)$ if and only if there is a distribution $T\in\OH$, the support of which has positive distance to all coordinate hyperplanes, such  that $L(S) = S\star T$ for all $S\in\Dip)$. Here $(S\star T)\vp = S_y(T_x\vp(xy)$ for all $\vp\in\Di$. The eigenvalues are $m_\alpha=T_x\Big(\frac{\sigma(x)}{x^{\alpha+\be}}\Big)$. \rm

Moreover, we show that every Hadamard operator on $\Dip$ maps $\CN$ to $\CN$, that is, defines a Hadamard operator on $\CN$. On the other hand not every Hadamard operator on $\CN$ can be extended to an operator on $\Dip$.

We follow the notation in \cite{Vhad,Vconv}. The class of Hadamard operators in $\CN$, that is, of continuous linear operators which admit all monomials as eigenvectors is denoted by $M(\Rn)$. The Hadamard operators are given by distributions $T\in\cE'(\Rn)$ by means of the formula $(M_T \vp)(y) = T_x\vp(xy)$. The dual $\cE'(\Rn)$ is an algebra with respect to $\star$-convolution given by the formula $(T\star S)\vp=T_x S_y \vp(xy)$ where $xy=(x_1y_1,\dots,x_dy_d)$. $T\in \cE'(\Rn)$ defines a $\star$-convolution operator $N_T:S\mapsto S\star T$ and $N_T=M_T^*$, that is, the dual operator of $M_T$.

Differential operators of the form $P(\theta)$ where $P$ is a polynomial and $\theta_j=x_j \partial_j$ or, equivalently, of the form $\sum_\alpha c_\alpha x^\alpha \partial^\alpha$ are called Euler operators and $\theta_j$ is called an Euler derivative. On $C^\infty$ these are the Hadamard operators $M_T$ with $\supp T =\{\be\}$
where $\be=(1,\dots,1)$.

We use standard notation of Functional Analysis, in  particular, of distribution theory. For unexplained notation we refer to \cite{DK}, \cite{LSI}, \cite{LS}, \cite{MV}.

The author thanks Pawe{\l} Doma\'nski, who passed away too early, for useful discussions on the topic of this paper.

\section{Basic properties}

\begin{definition}\label{d1} A map $L\in L(\Dip)$ is called a Hadamard operator if it admits all monomials as eigenvectors. The set of Hadamard operators we denote by $\cM(\Rn)$.
\end{definition}

Since the condition means that $L(x^\alpha)\in \span\{x^\alpha\}$ for all $\alpha\in\N_0^d$ the set $\cM(\Rn)$ is a closed subalgebra in $L_\sigma(\Dip)$ and therefore also in $L_b(\Dip)$. Here $\sigma$ denotes the topology of pointwise convergence and $b$ the topology of uniform convergence on bounded sets.

We define $m_\alpha$ by $L(x^\alpha)= m_\alpha x^\alpha$. Since the polynomials are dense in $\Dip$ the operator $L\in\cM(\Rn)$ is uniquely determined by the family $m_\alpha$, $\alpha\in\N_0^d$, of eigenvalues. Clearly the set $\Lambda(\Rn)$ of eigenvalue families is an algebra and $L\to m_\alpha$ is an algebra isomorphism.

To study and characterize the Hadamard operators in $\Dip$
we first need some preparations. We set $\R_*:=\R\setminus \{0\}$. For $a\in\NZ$ we define the dilation operator $D_a\in L(\cD'(\Rn))$ by
$$(D_a T)\vp:= T_x\left(\frac{\sigma(a)}{a_1\dots a_d}\,\vp\Big(\frac{x}{a}\Big )\right)$$
for $T\in\cD'(\Rn)$ and $\vp\in\cD(\Rn)$. By direct verification we see that $D_a\xi^\alpha=a^\alpha\xi^\alpha$.

\begin{lemma}\label{l1} For $L\in\cM(\Rn)$ and $a\in\R_*^d$ we have $L\circ D_a=D_a\circ L$.
\end{lemma}

\begin{proof} For any $\alpha$ we have $(L\circ D_a)\xi^\alpha=a^\alpha m_\alpha \xi^\alpha = (D_a\circ L)\xi^\alpha$. So the claim is shown for all polynomials and these are dense in $\cD'(\Rn)$. \end{proof}

\begin{notation} If $L\circ D_a=D_a\circ L$ for all $a$ we say: $L$ commutes with dilations. \end{notation}

By definition of $D_a$ we obtain for the dual map $D_a^*\in L(\cD(\Rn))$ of $D_a$:
\begin{lemma}\label{l2} For $a\in\R_*^d$ and $\vp\in\cD(\Rn)$ we have
$$(D_a^*\vp)(\xi) =\frac{\sigma(a)}{a_1\cdots a_d}\,\vp\left(\frac{\xi}{a}\right).$$
\end{lemma}

If $L$ commutes with dilations then $D_a^*\circ L^*= L^*\circ D_a^*$ for all $a\in\R_*^d$. For $\vp\in\cD(\Rn)$ we set $\psi = L^*\vp$ and obtain
$$\frac{\sigma(a)}{a_1\cdots a_d}\,\psi\left(\frac{x}{a}\right)=L^*_\xi\left(\frac{\sigma(a)}{a_1\cdots a_d}\,\vp\left(\frac{\xi}{a}\right)\right)[x].$$
For $\eta\in\R_*^d$ we set $a=1/\eta$ and obtain $\psi(\eta x)=L_\xi^*(\vp(\eta\xi))[x]$.

We have shown:
\begin{lemma}\label{l3} If  $M=L^*\in L(\cD(\Rn))$ and $L$ commutes with dilations, then
$$M_\xi(\vp(\eta \xi))[x]=(M\vp)(\eta x)$$
for all $\vp\in\cD(\Rn)$ and $\eta\in\R_*^d$.
\end{lemma}

For $\vp\in\cD(\Rn)$ we define now
$$T\vp= (M\vp)(\be)=(L\delta_\be)(\vp).$$
Then $T\in\cD'(\Rn)$ and for all $\eta\in\R_*^d$ we have
\begin{equation}\label{e1}(M\vp)(\eta)=T_\xi \vp(\eta\xi).\end{equation}

The problem with the right hand side of equation (\ref{e1}) is that for $\eta_j=0$, in general, $T$ cannot be applied to the non compact support function $\xi\to \vp(\eta\xi)$. For $T$ as above, however, the function $\eta\to T_\xi \vp(\eta\xi)$, $\xi\in\NZ$, is the restriction of a function in $\Di$.

\begin{definition}\label{d2} By $\DH$ we denote the set of distributions $T\in\Dip$ such that for every $\vp\in\Di$ the function
$y\to T_\xi \vp(\xi y)$, $y\in\NZ$, is the restriction of a function in $\Di$.
\end{definition}

This means that $T$ must have the following properties:
\begin{itemize}
\item[(*)] For every $\vp\in\Di$ there is $r>0$ such that $T_\xi \vp(\xi y)=0$ for $y\in\NZ$ with $|y|_\infty>r$,
\item[(**)] For every $\vp\in\Di$ the map $y\mapsto T_\xi \vp(\xi y)$, $y\in\NZ$, extends to a function in $\CN$.
\end{itemize}

For $T\in\DH$ we denote by $M_T$ the map which assigns to $\vp\in\Di$ the continuous extension of $y\mapsto T_\xi \vp(\xi y)$. From the closed graph theorem it follows easily:
\begin{lemma}\label{l8} $M_T\in L(\Di)$ for every $T\in\DH$.
\end{lemma}

We obtain the following representation theorem:

\begin{theorem}\label{t1} $L\in L(\Dip)$ commutes with dilations if and only if there is $T\in\DH$ such that $L(S)=S\star T$ for all $S\in\Dip$, here $(S\star T)\vp = S_y(T_x \vp(xy))$ for $\vp\in\Di$. In this case $T=L(\delta_\be)$.
\end{theorem}

\begin{proof} The assertion can be written as $L=M_T^*$. If $L$ commutes with translations then, by the above, $L^*=M_T$ for $T=L(\delta_\be)$ and $T\in\DH$ by formula (\ref{e1}). This formula also implies the result. \end{proof}

\begin{notation} For $T\in\DH$ we set $L_T(S)= S\star T$ for all $S\in\Dip$. \end{notation}

\begin{corollary}\label{c1} If $L\in\cM(\Rn)$ then there is $T\in\DH$ such that $L=L_T$.
\end{corollary}

\section{Properties of $\DH$}

First we will exploit the fact that $M_T\in L(\cD(\Rn))$.
For $\eps>0$ we set $$W_\eps=\{x\in\Rn\,:\, \min_j |x_j|\ge\eps\}$$
and we will use the following notations:\\
  For $r=(r_1,\dots,r_d)$, where all $r_j\ge 0$, we set $B_r=\{x\in\R^d:|x_j|\le r_j\text{ for all }j\}$ and for $s\ge 0$ we set $B_s:=B_{s\be}=\{x\in\Rn\,:\,|x|_\infty\le s\}$.\\
For $r=(r_1,\dots,r_d)\in\Rn$ and $s\in\R$ we set $r+s=r+s\be=(r_1+s,\dots,r_d+s)$.

\begin{lemma}\label{l4} 1. If $T\in\DH$ then there is $\eps>0$ such that $\supp T\subset W_\eps$.\\
2. If $T\in\Dip$ and there is $\eps>0$ such that $\supp T\subset W_\eps$ then $T$ satisfies \rm{(*)}.
\end{lemma}

\begin{proof}  1. If $T\in\DH$ there is $r>0$ such that $\supp M_T \vp\subset B_r$ for any $\vp\in\cD(B_1)$ and this implies, that $T_\xi\vp(\eta \xi)=0$ for any $\vp\subset \cD(B_1)$, $\eta\in\NZ$ and $|\eta|_\infty>r$.

We set $\eps=1/r$ and assume that $\supp \vp\cap W_\eps=\emptyset$. That is,
$$\supp \vp \subset\bigcup_j\, \{x\in\Rn\,:\, |x_j|<\eps\}.$$
Then we can write $\vp=\sum_j \vp_j$ with $\vp_j\in\cD(\{x\in\Rn\,:\, |x_j|<\eps\})$.

We fix $j$ and choose $\eta\in\R_*^d$ such that $\sup\{|x_\nu\eta_\nu|\,:\,x\in\supp\vp_j\}=1$ for all $\nu$. We set $\psi(\xi)=\vp_j(\xi/\eta)$. Then $\supp \psi\subset B_1$ and $|\eta|_\infty>r$ and therefore we have $T\vp_j= T_\xi\psi(\eta \xi)=0$.

Since this holds for every $j$, the proof of 1. is complete, the proof of 2. is obvious. \end{proof}

A special case is that of distributions with compact support.
\begin{corollary}\label{c2} $\cE'(\Rn)\cap\DH= \cE'(\NZ)$.
\end{corollary}

\begin{proof} This follows from Lemma \ref{l4} and the fact that (**) is fulfilled for distributions with compact support. \end{proof}

This will be used in Section \ref{s3} to handle the case of $T$ with compact support.

After having settled property (*) of $T\in\DH$ we turn to property (**). It is quite restrictive.

\begin{lemma}\label{l7} $\DH\subset \cS'(\Rn)$.
\end{lemma}

\begin{proof} Let $T\in\DH$, we may assume that $\supp T\subset W_2$. For $k\in\N_0^d$ we set $|\vp|_k = \|\vp^{(k)}\|_{L_1}$ and remark that for every $\cD(B_r)$ these norms are a fundamental system of seminorms.
On $\cD(B_\be)$ the family of distributions $T_{(y)}\vp:=T_x \vp(xy)$, $y\in\NZ\cap B_\be$, is weakly bounded hence equicontinuous. This means that there is $k\in\N_0^d$ and $C>0$ such that
$$|T_{(y)}\vp|\le C \|\vp\|_k,\quad \vp\in \cD(B_\be), |y|_\infty\le 1.$$
For $r$ with $r_j\ge 1$ for all $j$ and $\vp\in\cD(B_r)$ we set $\psi(x)=\vp(rx)$. Then $\psi\in \cD(B_\be)$ and therefore
$$|T\vp| = |T_{(1/r)} \psi|\le C \|\psi\|_k = C r^k \|\vp\|_k.$$

For every $r$ there is $t_r\in L_\infty(B_r)$, $\|t_r\|_\infty\le C r^k$ such that $T\vp=\int t_r(x) \vp^{(k)}(x) dx$ for all $\vp\in\cD(B_r)$.

We restrict us now to $r\in\N^d$ and set $U_r=[r_1,r_1+1[\times \dots \times [r_d,r_d+1[$. Then $\bigcup_{r\in\N^d}U_r = W_1$. We choose $\chi\in\cD(B_{1/2})$ with $\int \chi =1$ and set $\chi_r(x)=\int_{U_r}\chi(x-\xi) d\xi$. We obtain:
\begin{eqnarray*}T\vp &=& \sum_r T(\chi_r \vp) =\sum_r \int t_{r+2}(x) (\chi_r \vp)^{(k)}(x) dx\\
&=& \sum_r \int t_{r+2}(x) \left(\sum_\nu {k \choose \nu} \chi_r^{(k-\nu)}(x) \vp^{(\nu)}(x)\right) dx\\
&=& \sum_\nu {k \choose \nu} \int \vp^{(\nu)}(x) \left(\sum_r t_{r+2}(x) \chi_r^{(k-\nu)}(x)\right) dx\\
&=& \sum_\nu {k \choose \nu}\int \vp^{(\nu)}(x) \tau_\nu(x) dx.
\end{eqnarray*}

We have to  estimate the functions $\tau_\nu$. We set $\gamma(x)=\{r\,:\,x\in \prod_j[r_j-\frac{1}{2},r_j+\frac{3}{2}[\}$
\begin{eqnarray*}
|\tau_\nu(x)| &\le&  \sum_{r\in\gamma(x)} |t_{r+2}(x)|\, |\chi_r^{(k-\nu)}(x)|.
\end{eqnarray*}
For all $2^d$ elements $r\in \gamma(x)$ we have $|t_{r+2}(x)|\le C (r+2)^k$ and $r\le |x|+\frac{1}{2}$ (here $|x|=(|x_1|,\dots,|x_d|)$) and therefore
$|t_{r+2}(x)|\le C(|x|+3)^k$. With a new constant $C_\nu$ we have
$$|\tau_\nu(x)|\le C_\nu |x|^k.$$
This shows the result. \end{proof}

The necessary conditions we have found are far from being sufficient as the following example shows

\begin{example}\label{ex1} Let $d=1$, we set $T\vp=\int_1^\infty \vp(x)dx$, then $T_x\vp(x y)= \frac{1}{y}\int_y^\infty \vp(x)dx$ for all $y>0$ which, in general, is unbounded near $0$.
The distribution is $T\in\cS'(\R)$ and has support in $W_1$.
\end{example}

The following definition is in analogy to a characterization of the space $\cO_C'$ of rapidly decreasing distributions in L. Schwartz \cite{LS}, \S 5, p. 100.

\begin{definition}\label{d3} $T\in \OH$ if for any $k$ there are finitely many functions $t_\beta$ such that $(1+|x|^2)^{k/2}t_\beta\in L_\infty(\Rn)$ and such that $T=\sum_\beta \theta^\beta t_\beta$.
\end{definition}

\begin{proposition}\label{l5} If $T\in\OH$ then for every $\vp\in\cD(\Rn)$ the function $y\mapsto T_x(\vp(xy))$, $y\in\R_*^d$, extends to a function in $C^\infty(\Rn)$.
\end{proposition}

\begin{proof} We have to show that for every $\alpha$ the function $y\mapsto \partial_y^\alpha T_x(\vp(xy))$, $y\in\R_*^d$, extends to a continuous function on $\Rn$.

By definition of $\OH$ we have to show our proposition only for $T=\theta^\beta t$ where $(1+|x|^2)^{k/2}t\in L_\infty(\Rn)$ and $k$ is suitably chosen. Since the formal adjoint $(\theta^*)^\beta$ of the Euler-operator $\theta^\beta$ is again an Euler-operator we have
$$\int (\theta^\beta t)(x) \vp(x) dx = \int t(x) ((\theta^*)^\beta \vp)(x) dx = \sum_\nu c_\nu \int t(x) x^\nu \vp^{(\nu)}(x) dx$$
where the sum is finite with suitable $c_\nu$. Therefore it is enough to study the function
\begin{equation}\label{e3}F(y):= \int t(x) (xy)^\nu \vp^{(\nu)}(xy) dx, \quad y\in\NZ.\end{equation}
We have to show that all limits $\lim_{y\to y_0} F^{(\alpha)}(y)$ , $y_0\in \Rn$, $\alpha\in\N_0^d$ exist.
\begin{multline}\label{e4}F^{(\alpha)}(y) = \int t(x) x^\alpha  (x^\nu \vp^{(\nu)}(x))^{(\alpha)}[xy] dx\\
= \sum_{\alpha-\nu\le\gamma\le\alpha} c_\gamma \int t(x) x^\alpha (xy)^{\nu-\alpha+\gamma} \vp^{(\nu+\gamma)}(xy) dx.
\end{multline}
If $k$ was chosen so large that $t(x) x^\alpha \in L_1(\Rn)$ then the limits exist as requested.  \end{proof}

Therefore we have shown:

\begin{theorem}\label{p1} If $T\in\OH$ and $\supp T\subset W_\eps$ for some $\eps>0$, then $T\in\DH$.
\end{theorem}

To show the inverse we give a description of $\DH$ in terms of the regularizations of a distribution.

\begin{theorem}\label{t8} $T$ satisfies (**) if and only if $\,T*\chi$ satisfies (**) for all $\chi\in L_1(B_1)$.\\
In this case the set of maps $\{\vp\mapsto (T*\chi)_x \vp(x\,\cdot)\,:\, \chi\in L_1(B_1),\,\|\chi\|_{L_1}\le 1\}$ is equicontinuous in $L(\Di,\CN)$.
\end{theorem}

\begin{proof} We assume that $T$ satisfies (**) and want to show that the same holds for $T*\chi$. We need some preparation.
We set:
$$\psi(\xi,y_2)=\int\chi(x_2) \vp(\xi+x_2 y_2) dx_2$$
and remark that for $\vp\in\cD(B_r)$ and $\chi\in L_1(B_1)$ we have $\psi(\cdot,y_2)\in\cD(B_{r+|y_2|})$.  This implies that the map which assigns to every $\vp\in\Di$ the function $y_2\mapsto \psi(\cdot,y_2)$ is a continuous linear map $\Phi:\Di\to C^\infty(\Rn,\Di)$.

We set
$$F(y_1,y_2) = T_{x_1} \psi(x_1 y_1,y_2) \in C^\infty(\R_*^d\times\R^d).$$
By assumption it extends to a function $\hat{F}$ on $\Rn\times\Rn$ such that $\hat{F}(\cdot,y_2)$ is in $\CN$ for every $y_2\in\Rn$.

The map which assigns to every $g\in\Di$ the continuous extension of $T_\xi g(\xi y_1)$ defines a continuous linear map $\Psi:\Di\to\CN$. We obtain
$$\hat{F}(y_1,y_2)=\Psi\{\Phi(\vp)[y_2]\}[y_1]\in C^\infty(\Rn,C^\infty(\Rn))=C^\infty(\Rn\times\Rn).$$
Therefore $\hat{F}(y,y)\in\CN$ and it is the extension of
$$F(y,y)=T_{x_1}\int\chi(x_2) \vp((x_1+x_2)y)dx_2= (T*\chi)_x\vp(xy).$$

This proves one direction of the Theorem, it remains to show the other implication. We will use the idea of proof in \cite[\S 7, Th\'eor\`eme XX]{LS}.

First we show the additional assertion of the theorem. We assume that $\,T*\chi$ satisfies (**) for all $\chi\in L_1(B_1)$. We consider the map $L_1(B_1)\to L(\Di,\CN)$ defined by $\chi\mapsto [\vp\mapsto (T*\chi)_x \vp(x\,\cdot)]$ (cf. Lemma \ref{l8}). If $\|\chi\|\to 0$ and $[\vp\mapsto (T*\chi)_x \vp(x\,\cdot)]\to A$ in $L(\Di,\CN)$  then for fixed $y\in\NZ$ and all $\vp$ we have $(T*\chi)_x \vp(xy)\to 0$ and therefore $(A\chi)\vp=0$ on $\NZ$, hence on $\Rn$. So the map $L_1(B_1)\to L(\Di,\CN)$ has closed graph and, by de Wilde's Theorem, is continuous. This shows the assertion.

We fix $\chi\in\Di$ with $\int |\chi| =1$. Then for every $\vp\in\Di$ the function $(T*\chi)_x \vp(xy)$, $y\in\NZ$, extends to a function in $\CN$ and
the set $\{F_\chi:\vp\mapsto (T*\chi)_x \vp(xy)\,:\,\chi \in \Di,\,\chi\ge 0, \int\chi=1\}$ is equicontinuous. Therefore it is relatively compact in $L(\cD(B_R), \CN)$. We fix $\chi$ and set $\chi_\eps(x)=\eps^{-d} \chi(x/\eps)$ for $\eps>0$. Then there is a sequence $\eps_n \downarrow 0$ such that $F_{\chi_{\eps_n}}$ converges to some  $F \in L(\cD(B_R),\CN)$. Since
$$F_{\chi_{\eps_n}}=T_\xi\Big(\int \chi_{\eps_n}(x)\vp((x+\xi)y)dx\Big) \to T_\xi \vp(\xi y)$$
for every $y\in\NZ$ we see, that $T_\xi \vp(\xi y)$ extends to a function in $\CN$. \end{proof}

Now we can show the inverse of Theorem \ref{p1}.

\begin{theorem}\label{l10}If $T\in\DH$ then for every $\beta>0$ there is a function $t_\beta$ such that $x^\beta t_\beta$ is bounded and an Euler operator $P(\theta)$ such that $T = P(\theta) t_\beta$. In particular $T\in\OH$.
\end{theorem}

\begin{proof}
Let $T\in \DH$ with $\supp T\subset W_{2\eps_0}$. Let $\chi\in\Di$ with $\chi\ge 0$, $\int \chi =1$ and $\supp \chi\subset B_{\eps_0}$. Then for every $\vp\in \cD(B_\be)$ the function $(T*\chi)_x \vp(xy)$, $y\in\NZ$, extends to a function in $\Di$. We set $\tau = T*\chi$.

We set for $y\in\NZ$ and $\vp\in\Di$:
$$F(y) = \int \tau(x) \vp(xy) dx.$$
Then $F\in C^\infty(\NZ)$, $F(y)=0$ for $|y|_\infty>1/\eps_0$. For $\beta\in\N_0$, we have
$$F^{(\beta)}(y) = \int \tau(x) x^\beta \vp^{(\beta)}(xy) dx.$$

By Theorem \ref{t8} on $\cD(B_\be)$ the set of distributions $\vp\to F^{(\beta)}(y)$, $y\in\NZ$,  $\chi$ as above, is equicontinuous. Hence there is $p$ such that all these distributions extend to $\cD^{|\beta|+p}(B_\be)$ and the set of these distributions is bounded in $\cD^{|\beta|+p}(B_\be)'$.

For $\alpha\in\N_0$ we choose $\vp_\alpha\in \cD^{\alpha}[0,+1]$ such that $\vp_\alpha\in C^\infty(\R\setminus\{0,1\})$  and
$$\vp_\alpha(x) = x^{\alpha+1}\,\frac{(x-1)^{\alpha+1}}{(2\alpha+2)!} \text{ for } 0\le x \le 1.$$
We define for $\alpha\in\N_0^d$ such that $\cD^{|\alpha|}(B_\be)\subset \cD^{|\beta|+p}(B_\be)$
$$\vp_\alpha(x)= \prod_{j=1}^d \vp_{\alpha_j}(x_j)$$
and consider $F(y)=\int \tau(x) \vp_\alpha(xy) dx.$
Then we have for $y\in (0,+\infty)^d$, setting $\alpha+1=(\alpha_1+1,\dots,\alpha_d+1)$, etc.
$$F^{(2\alpha+2)}(y) = \int_{[0,1/y]^d} \tau(x) x^{2\alpha+2} dx$$
and therefore
$$F^{(2\alpha+3)}(y) = -\frac{1}{y^{2\alpha+4}}\, \tau\Big(\frac{1}{y}\Big).$$

With an analogous argument for the other `{quadrants}' we get for general $y\in\NZ$
$$F^{(2\alpha+3)}(y) = -\frac{\sigma(y)}{y^{2\alpha+4}}\, \tau\Big(\frac{1}{y}\Big).$$

We set $G(x)=\sigma(x)F^{(\beta)}(\frac{1}{x})$. Then $G$ is a bounded function on $\NZ$ with a bound independent of $\chi$ and $\sigma(x)\,F^{(\beta)}(x)=G\Big(\frac{1}{x}\Big)$. We calculate the derivatives. With certain coefficients $c_\nu$ we have:
$$\sigma(x)\,F^{(\beta+q)}(x) = {\sum_\nu}' c_\nu \, \frac{1}{x^{\nu+q}}\, G^{(\nu)}\Big(\frac{1}{x}\Big)$$
where ${\sum}'_\nu$ runs over $1\le\nu_j\le q$. We choose $q=2\alpha + 3 -\beta$. Replacing $x$ with $1/x$ and $c_\nu$ with $-c_\nu$ we obtain:
\begin{equation}\label{eq1}
x^{\beta+1}\, \tau(x) = {\sum_\nu}'c_\nu\, x^\nu\, G^{(\nu)}(x).
\end{equation}

We analyze now functions of type $x^{m+\nu} G^{(\nu)}(x)$ and claim

\begin{lemma}\label{l9} Any such function can be written as a linear combination of the functions $(x^{m+j} G(x))^{(j)}$, $0\le j_i\le\nu_i$ for all $i$.
\end{lemma}

\begin{proof} For any $m\in\Z$ there are constants $c_{m,\nu,j}$ such that
$$(x^{m+\nu} G(x))^{(\nu)} = \sum_{j=0}^\nu c_{m,\nu,j}\, y^{m+\nu-j} G^{(\nu-j)}(x)$$
and therefore,since $c_{m,\nu,0}=1$,
$$x^{m+\nu} G^{(\nu)}(x)= (x^{m+\nu} G(x))^{(\nu)} -  {\sum_j}' c_{m,\nu,j}\, x^{m+\nu-j} G^{(\nu-j)}(x).$$
where ${\sum_j}'$ runs over all $j$ with $j\neq 0$ and $0\le j_i\le \nu_i$ for all $i$. Induction over $|\nu|$ yields the result. \end{proof}

{\bf End of the proof of Theorem \ref{l10}:} From Lemma \ref{l9} and equation (\ref{eq1}) we obtain with new constants $c_j$:
$$\tau(x)=\sum_{|j|_\infty\le p+2} c_j\, (x^j\,(x^{-\beta-1}G(x)))^{(j)}.$$

We fix now $\chi\in\Di$ with $\chi\ge 0$ and $\int\chi =1$. For $\eps>0$ we set $\chi_\eps=\eps^{-d} \chi(x/\eps)$, $\tau_\eps = T*\chi_\eps$ and $G_\eps$ the corresponding function. Then we have
\begin{equation}\label{eq4}(T*\chi_\eps)\vp = \int \tau_\eps(x)\vp(x) dx = \sum_{|j|_\infty\le p+2} c_j \int x^{-\beta-1} G_\eps(x) x^j \vp^{(j)}(x) dx.\end{equation}
We have $\lim_{\eps\to 0} (T*\chi_\eps)\vp= T\vp$. On the other hand $\{G_\eps\,:\,\eps>0\}$ is bounded in $L_\infty(\Rn)=L_1(\Rn)'$. Therefore there is $G\in L_\infty(\Rn)$ and a sequence $G_{\eps_n}$ which converges to $G$ in the weak$^*$-topology with respect to $L_1(\Rn)$.

From equation (\ref{eq4}) then follows:
\begin{eqnarray*}T\vp &=& \sum_{|j|_\infty\le p+2} c_j \int x^{-\beta-1} G(x) x^j \vp^{(j)}(x) dx = \int (x^{-\beta-1} G(x)) (P(\theta) \vp)(x)dx\\
&=& \int(P(\theta)^*(x^{-\beta-1} G(x)) \vp(x) dx.
\end{eqnarray*}

 Putting $t_\beta:= x^{-\beta-1} G(x)$ we have shown: For every $\beta$ there is a function $t_\beta$  such that $x^\beta t_\beta\in L_\infty(\Rn)$ and an Euler operator $Q(\theta):=P(\theta)^*$ such that $T=Q(\theta)t_\beta$.
This completes the proof. \end{proof}

From Theorems \ref{p1} and \ref{l10} we obtain one of the main results of this paper:
\begin{theorem} \label{t11} $\DH=\{T\in\OH\,:\, \supp T\subset W_\eps \text{ for some }\eps>0\}$.
\end{theorem}

\section{The space $\OH$}

We recall the definition of the space $\OH$:

$T\in \OH$ if for any $k$ there are finitely many functions $t_\beta$ such that $(1+|x|^2)^{k/2}t_\beta\in L_\infty(\Rn)$ and such that $T=\sum_\beta \theta^\beta t_\beta$.

The space $\OC$ of L. Schwartz may be defined by any of the following equivalent properties (see \cite[\S 5, Th\'eor\`eme IX]{LS}). Let $T\in\Dip$, then $T\in\OC$ if and only if 1. or 2.:
\begin{enumerate}
\item For any $k$ there are finitely many functions $t_\beta$ such that $(1+|x|^2)^{k/2}t_\beta\in L_\infty(\Rn)$ and such that $T=\sum_\beta \partial^\beta t_\beta$.
\item For any $\chi\in\Di$, $T*\chi$ is a rapidly decreasing continuous function.
\end{enumerate}

Without proof we admit.
\begin{lemma}\label{l11} $\cE'(\Rn)\subset \OC\cap\OH$.
\end{lemma}
And this implies:
\begin{lemma}\label{l12} If $T|_{\Rn\setminus B_R}=S|_{\Rn\setminus B_R}$ for some $R>0$ and $S\in\OH$ or $S\in\OC$ then $T\in\OH$ or $T\in\OC$, respectively.
\end{lemma}

\begin{proof} $T=S+(T-S)$ and $T-S\in\cE'(\Rn)$. \end{proof}

\begin{proposition}\label{p4} 1. $\Oc\subset\Oh$.\\
2. If $\supp T\subset W_\eps$ for some $\eps>0$ and $T\in\OC$ then $T\in\OH$.
\end{proposition}

\begin{proof} 2. It is enough to show the claim for $T=t_\beta^{(\beta)}$, $(1+|x|^2)^{k/2} t_\beta$ bounded. Set $\tau_\beta = \frac{1}{x^\beta} t_\beta$. Then $\tau_\beta$ is decreasing even faster and $(x^\beta \tau_\beta)^{(\beta)}=t_\beta^{(\beta)}$. Since $f\mapsto (x^\beta f)^{(\beta)}$ is an Euler operator the proof of 2. is complete.\\
1. Choose $\chi\in\cD[-1,+1]$, $\chi\equiv 1$ in a neighborhood of $0$. For $T\in\Oc$ set $S=(1-\chi)T$. Then $\supp S\subset W_\eps$ for some $\eps>0$ and, due to Lemma \ref{l12}, $S\in\Oc$. By part 2. of this Proposition we have $S\in\Oh$ and therefore, again by Lemma \ref{l12}, $T\in\Oh$. \end{proof}

The space $\Oc$ is a proper subspace of $\Oh$, as the following example shows.

\begin{example}\label{ex4} If $T=e^{-ix}$, that is, $T\vp = \int e^{-ix} \vp(x) dx$, then\\
1. $T\not\in \Oc$,\,
2. $T\in \Oh$.
\end{example}

\begin{proof} 1. Let $\chi\in\cD(\R)$, then $(T*\chi)(x)=\int e^{-i\xi} \chi(x-\xi) d\xi= \hat{\chi}(-1)\, e^{-ix}$. By the second definition of $\Oc$ (see above) $T\not\in\Oc$.\\
2. To show that $T\in\Oh$ we choose $\chi\in\cD[-1,+1]$, $\chi\equiv 1$ in a neighborhood of $0$. For given $k$ we set $t_k(x)=\frac{i^k}{x^k} (1-\chi(x))e^{-ix}$. Then $(1+x^2)^{k/2} t_k(x)$ is bounded and
$(x^k t_k(x))^{(k)}=(i^k(1-\chi(x)) e^{-ix})^{(k)}= e^{-ix} + g(x)$ where $g$ has compact support.
Hence $T= (x^k t_k(x))^{(k)} - g$ where $g$ has compact support. This shows the result like above. \end{proof}

By Proposition \ref{l5} we know now that for $T$ as in Example \ref{ex4} the function $T_x\vp(xy)$, $y\in \R_*$, extends to a $C^\infty$-function on $\R$. For this example we can make it explicit, even for higher dimensions, setting $e^{-ix}=e^{-i(x_1+\dots+x_d)}$. For $\vp\in\Di$ and $y\in \NZ$ we set
$F(y)=T_x\vp(xy)$. We obtain for $y\in\NZ$
\begin{eqnarray*}
\partial^\alpha F(y)&=& \int e^{-ix} x^\alpha \vp^{(\alpha)}(xy) dx =\frac{\sigma(y)}{y^{\alpha+\be}}\int e^{-i(x/y)} x^\alpha \vp^{(\alpha)}(x) dx\\ &=& \frac{\sigma(y)}{y^{\alpha+\be}}\,  \widehat{x^\alpha \vp^{(\alpha)}}\Big(\frac{1}{y}\Big).
\end{eqnarray*}
Since $\widehat{x^\alpha \vp^{(\alpha)}}\in\cS(\Rn)$ we obtain $\lim_{y\to y_0} F^{(\alpha)}(y)=0$ for every $y_0\in\Rn\setminus \NZ$. That means, if we denote the extended function again by $F$, that $F^{(\alpha)}(y)=0$ on all coordinate hyperplanes and for all $\alpha$.

Returning to the one-dimensional case we present another example which we take from \cite[\S 5, p. 100, (VII,5;1)]{LS}.
\begin{example}\label{ex5} If $T=e^{i\pi x^2}$ then $T\in \Oc$ and therefore $T\in\Oh$. The function $e^{i\pi x^2}$ is bounded, but its derivatives are not.
\end{example}

\section{Eigenvalues}\label{s3}

In this section we study Hadamard operators in terms of the representing distribution in $\DH$. A special case are the distributions in $\DH$ with compact support.

In consequence of Lemma \ref{l1}, Theorem \ref{t1} and Corollary \ref{c2} we obtain:

\begin{theorem}\label{t3} For $T\in\Dip$ the following are equivalent:\\
1. $T\in\cE'(\Rn)$ and the $\star$-homomorphism $N_T$ can be extended to a map in $\cM(\Rn)$.\\
2. $T\in \DH$ and $L_T\cE'(\Rn)\subset\cE'(\Rn)$.\\
3. $T\in \cE'(\R_*^d)$.

In this case $m_\alpha=T_x\Big(\frac{\sigma(x)}{x^{\alpha+\be}}\Big).$
\end{theorem}

\begin{proof} 1. $\Rightarrow$ 2.: By assumption there is $L\in \cM(\Rn)$ such that $L(S)=N_T(S)=S\star T$ for all $S\in\cE'(\Rn)$. Then $T=N_T(\delta_\be)=L(\delta_\be)\in\DH$. By definition $N_T(S)=S\star T = L_T(S)$ for $S\in \cE'(\Rn)$.

2. $\Rightarrow$ 3.: By assumption $T=L_T(\delta_\be))\in\cE'(\Rn)$. Hence $T\in \cE'(\Rn)\cap\DH= \cE'(\NZ)$, by Corollary \ref{c2}.

3. $\Rightarrow$ 1. By Corollary \ref{c2} $T\in\cE'(\R_*^d)\subset \DH$ and we have $L_T(S)=S\star T=N_T(S)\in \cE'(\Rn)$ for $S\in\cE'(\Rn)$.

Moreover
\begin{eqnarray*}L_T(\xi^\alpha)[\vp]&=&\int \xi^\alpha ( T_x \vp(x\xi)) d\xi=T_x\Big(\frac{1}{x^{\alpha}}\int (x\xi)^\alpha \vp(x\xi)d\xi\Big)\\
&=& T_x\Big(\frac{\sigma(x)}{x^{\alpha+\be}}\int\eta^\alpha \vp(\eta) d\eta\Big)\\
&=& \int (m_\alpha \eta^\alpha) \vp(\eta)d\eta.
\end{eqnarray*}
This shows that $L_T(\xi^\alpha)=m_\alpha \xi^\alpha$ with $m_\alpha=T_x\Big(\frac{\sigma(x)}{x^{\alpha+\be}}\Big).$ \end{proof}

Notice that the interchange of $T$ and the integral is clear since $T\in\cE'(\Rn)$. For $T$ with non-compact support this is more complicated. We further study operators represented by distributions in $\DH$. By Theorems \ref{t1} and \ref{t11} we know that they define operators $L_T\in L(\Dip)$ commuting with dilations. We show that these are, in fact, Hadamard operators.

\begin{theorem}\label{t7} If $T\in\DH$, then $L_T$ is a Hadamard operator with eigenvalues $m_\alpha=T_x\Big(\frac{\sigma(x)}{x^{\alpha+\be}}\Big).$
\end{theorem}

\begin{proof} It remains to show that $L_T$ admits all monomials as eigenvectors. It is sufficient to assume that $T = (-1)^{|k|} (x^k \tau)^{(k)}$, where $\tau\in L_1(\Rn)$ and $k\in\N_0^d$. Then the function
$$f(x,\xi)= \xi^\alpha \tau(x) \xi^k \vp^{(k)}(x\xi)$$
is in $L_1(\Rn\times\Rn)$. To see this let $\supp \vp \subset B_R$. Then $f(x,\xi)\neq 0$ only if $|x_j|\ge\eps$ and $|x_j \xi_j|\le R$ for all $j$, hence only for $|\xi_j|\le R/\eps$ for all $j$.

\begin{eqnarray*} \int\xi^\alpha T_x\vp(x\xi) d\xi &=& \int \xi^\alpha\Big(\int \tau(x) (x\xi)^k \vp^{(k)}(x\xi) dx\Big) d\xi\\
&=& \int \tau(x) \frac{1}{x^\alpha} \Big(\int (x\xi)^{\alpha+k} \vp^{(k)}(x\xi) d\xi\Big)dx\\
&=& \int \tau(x) \frac{\sigma(x)}{x^{\alpha+\be}} \Big(\int \eta^{\alpha+k}\vp^{(k)}(\eta) d\eta\Big) dx\\
&=& \int \tau(x) \frac{\sigma(x)}{x^{\alpha+\be}} (-1)^{|k|} \frac{(\alpha+k)!}{\alpha!} \Big(\int \eta^\alpha \vp(\eta)d\eta\Big)dx.
\end{eqnarray*}
Therefore
$$\int\xi^\alpha T_x\vp(x\xi) d\xi=\int(m_\alpha \eta^\alpha) \vp(\eta) d\eta$$
where
\begin{eqnarray*}m_\alpha &=& (-1)^{|k|} \frac{(\alpha+k)!}{\alpha!} \int \tau(x) \frac{\sigma(x)}{x^{\alpha+\be}} dx.\\
\end{eqnarray*}
Taking into account that for $x$ with $\min_j |x_j|>0$
$$x^k \Big(\frac{\sigma(x)}{x^{\alpha+\be}}\Big)^{(k)}= (-1)^{|k|}\, \frac{(\alpha+k)!}{\alpha!}\, \frac{\sigma(x)}{x^{\alpha+\be}}$$
we finally obtain
$$m_\alpha= T\Big(\frac{\sigma(x)}{x^{\alpha+\be}}\Big)$$
which proves the result. \end{proof}

\sc Remark: \rm We needed in the proof only very weak assumptions on $T$. Under these $\vp\mapsto T_x \vp(xy)$ might not send $\Di$ into $\Di$ hence its adjoint sends something much smaller to $\Dip$, but it sends $\xi^\alpha$ to $m_\alpha \xi^\alpha$ for all $\alpha\in\N_0^d$.

\section{Hadamard operators in $\Dip$ and in $\CN$}

We study now the problem when a Hadamard operator $M$ on $\CN$ extends to an operator on $\Dip$ and, on the other side, when an operator $L\in\cM(\Rn)$ leaves $\CN$ invariant, that is, $L\CN\subset \CN$.

We start with the latter question, which has a rather straightforward answer.

\begin{theorem}\label{t4} If $T\in\DH$
then $L_T\in\cM(\Rn)$ and $L_T(\CN)\subset \CN$.
\end{theorem}

 \begin{proof} We may assume that $T=P(\theta) t$ where $t\in L_1(\Rn)$ with $\supp t\subset W_\eps$ for some $\eps>0$ and $P(\theta)$ is an Euler operator. Let $P^*(\theta)$ denote the formal adjoint of $P(\theta)$ which is again an Euler operator. Then  for $f\in\CN$, $\vp\in\Di$ the function
$$f(y)t(x)P^*(\theta)_x\vp(xy),\quad x,y\in\Rn$$
is in $L_1(\Rn\times\Rn)$. We obtain,using that Euler operators commute with dilations:
\begin{eqnarray*}
(L_T f)\vp &=& \int f(y)\Big(\int t(x) P^*(\theta)_x \vp(xy) dx\Big) dy\\
&=& \int\int  f(y) t(x) (P^*(\theta)\vp)(xy) dx dy.
\end{eqnarray*}
We apply the substitution $xy=\eta,\,y=\xi\eta$, its Jacobian determinant is $\frac{(-1)^d}{\xi_1\cdots\xi_d}$.
$$(L_T f)\vp = \int\Big(\int f(\xi\eta)\, t\Big(\frac{1}{\xi}\Big)\,\frac{\sigma(\xi)}{\xi_1\cdots\xi_d}d\xi\Big)(P^*(\theta)\vp)(\eta) d\eta.$$
We have shown that on $\CN$ we have $L_T=P(\theta)\circ M_{T^\#}$ where
$$T^\# = t\Big(\frac{1}{\xi}\Big)\,\frac{\sigma(\xi)}{\xi_1\cdots\xi_d}.$$
Notice that $t\Big(\frac{1}{\xi}\Big)\,\frac{\sigma(\xi)}{\xi_1\cdots\xi_d}$ is an $L_1(\Rn)$-function with compact support, hence $T^\#\in\cE'(\Rn)$. \end{proof}

\sc Remark: \rm Like in Theorem \ref{t7} we needed much weaker assumptions than $T\in\OH$, cf. the remark at the end of Section \ref{s3}. This corresponds to the fact that not all Hadamard operators in $\CN$ extend to operators in $\cM(\Rn)$ as the following proposition shows:

\begin{proposition}\label{ex2} If $\supp T=\{0\}$ and $T\neq 0$ then $M_T$ cannot be extended to a map in $\cM(\Rn)$
\end{proposition}

\begin{proof} If $T=\sum_\beta c_\beta \delta^{(\beta)}$ then $M_T(f)[x] = \sum_\beta c_\beta (-1)^{|\beta|} f^{(\beta)}(0)x^\beta$, hence $R(M_T)\subset E=\span\{x^\beta\,:\, \beta\in e\}$ where $e$ is a finite set. Assume $L\in\cM(\Rn)$ and $L|_{\CN}=M_T$. Since $E$ is closed in $\Dip$ and $\CN$ is dense in $\Dip$ we have $R(L)\subset E$ and this implies that $L(S)=\sum_\beta S(\vp_\beta)x^\beta$ for all $S\in\Dip$ where $\vp_\beta\in\Di$ for all $\beta\in e$. This implies $c_\beta (-1)^\beta f^{(\beta)}(0) = \int f(\xi) \vp_\beta(\xi) d\xi$ for all $\beta\in e$ and $f\in\CN$ which is possible only if $c_\beta=0$ for all $\beta\in e$, that is, $T=0$. \end{proof}

To study the problem which Hadamard operators on $\CN$ extend to operators in $\cM(\Rn)$ we consider an operator $M=M_T\in M(\Rn)$, $T\in\cE'(\Rn)$. We may assume that $T=(-1)^{|\beta|} t^{(\beta)}$ where $t\in L_1(\Rn)$ with compact support.

Here we assume that $f\in\CN$ and $\vp\in\Di$. Then the function
$$t(x)\vp(y) y^\beta f^{(\beta)}(xy),\quad x,y\in\Rn$$
is in $L_1(\Rn\times\Rn)$. We use the same coordinate transformation as above. Then the function
$$\frac{\sigma(\xi)}{\xi_1\cdots\xi_d}\, t\Big(\frac{1}{\xi}\Big) \vp(\xi\eta) (\xi\eta)^\beta f^{(\beta)}(\eta),\quad \xi,\eta\in\Rn$$
is again in $L_1(\Rn\times\Rn)$. By Fubini's theorem and the change of variables we obtain:
\begin{eqnarray*}
\int \vp(y)\left\{\int t(x) y^\beta f^{(\beta)}(xy) dx\right\} dy &=& \int\int t(x)\vp(y) y^\beta f^{(\beta)}(xy) dx dy\\
&=& \int\int \frac{\sigma(\xi)}{\xi_1\cdots\xi_d}\, t\Big(\frac{1}{\xi}\Big) \vp(\xi\eta) (\xi\eta)^\beta f^{(\beta)}(\eta) d\xi d\eta\\
&=& \int\left(\int \frac{\sigma(\xi)}{\xi_1\cdots\xi_d}\, t\Big(\frac{1}{\xi}\Big)\xi^\beta \vp(\xi\eta)d\xi \right) (P(\theta) f)(\eta)  d\eta\\
&=& \int (P(\theta)f)(\eta) S_\xi \vp(\xi\eta) d\eta
\end{eqnarray*}
where $P(\theta)f(\eta) = \xi^\beta f^{(\beta)}(\eta)$ and the distribution $S$ is defined by
$$S=\frac{\sigma(\xi)}{\xi_1\cdots\xi_d}\,\xi^\beta\,t\Big(\frac{1}{\xi}\Big).$$

Our problem now comes down to the question: when is $S\in\DH$? Then we have
$$M_T=L_S\circ f(\theta).$$
Since clearly $\supp S\subset W_\eps$ for some $\eps>0$, the question is: when does $\eta\mapsto S_\xi \vp(\xi\eta),\,\eta\in\NZ$ extend to a function in $\CN$?

Sufficient for that is if $ \frac{\xi^\gamma}{\xi_1\cdots\xi_d}\, t\Big(\frac{1}{\xi}\Big)\in L_1(\Rn)$ for all $\gamma\in\N_0^d$ and this is equivalent to
$$x^{-\gamma} t(x) \in L_1(\Rn)\text{ for all } \gamma\in\N_0^d.$$
We have shown:

\begin{theorem}\label{t5} If $T\in\cE'(\Rn)$ and there are finitely many functions $t_\beta\in L_1(\Rn)$ with compact support such that $T=\sum_\beta t_\beta^{(\beta)}$ and $x^{-\gamma} t_\beta(x) \in L_1(\Rn)$ for all $\gamma\in\N_0^d$ then $M_T:\CN\to\CN$ extends to a map in $\cM(\Rn)$.
\end{theorem}

A condition about the behaviour of $T$ at $0$ is necessary as we have seen on Proposition \ref{ex2}.

That in the assumptions for Theorem \ref{t5} we need a strong vanishing condition at zero and that a condition in the spirit of the $\OH$-condition is not sufficient is shown by the following easy example.

\begin{example}\label{ex3} For any $p\in\N_0$ we define a function $t_p$ as follows: $t_p(x)=0$ for $x<0$, $t(x)=\chi(x) x^p/p!$ for $x\ge 0$ where $\chi\in\Di$ and $\chi\equiv 1$ in a neighborhood of $0$. Then $\delta = t^{(p)}_p + t_0$ where $t_0\in \cD(]0,+\infty[)$. So for any $P$ the distribution $\delta$ is a sum of derivatives of functions with zeroes of order $p$ in $0$, but $M_\delta$ does not extend to an operator in $\cM(\Rn)$ (see Proposition \ref{ex2}).
\end{example}

\end{document}